\documentclass{article}%
\usepackage{amssymb}
\usepackage{amsmath}
\usepackage{amsfonts}
\usepackage{graphicx}%
\providecommand{\U}[1]{\protect\rule{.1in}{.1in}}
\newtheorem{theorem}{Theorem}
\newtheorem{acknowledgement}[theorem]{Acknowledgement}

\newtheorem{corollary}[theorem]{Corollary}

\newtheorem{lemma}[theorem]{Lemma}

\newtheorem{proposition}[theorem]{Proposition}
\newtheorem{remark}[theorem]{Remark}

\begin{document}

\begin{center}
{\LARGE {\Large On the quenched CLT for stationary Markov chains}.}

\bigskip

{\large Dedicated to Michael Lin's 80th birthday.}

\bigskip

Magda Peligrad

\bigskip

Department of Mathematical Sciences, University of Cincinnati, PO Box 210025,
Cincinnati, Oh 45221-0025, USA. \texttt{ }
\end{center}

email: peligrm@ucmail.uc.edu

\bigskip

\noindent\textit{Keywords:} Markov chains, quenched limit theorems, central
limit theorem, projective criteria.

\smallskip

\noindent\textit{Mathematical Subject Classification} (2010): 60F05, 60F15,
60J05, 60G10.

\bigskip

\textbf{Abstract}. In this paper we give sufficient conditions for the almost
sure central limit theorem started at a point, known under the name of
quenched central limit theorem. This is achieved by using a new idea of
conditioning with respect to both the past and the future of the Markov chain.
As applications we provide a new sufficient projective conditions for the
quenched CLT.

\section{Introduction and the main result}

We assume that $(\xi_{n})_{n\in\mathbb{Z}}$ is a stationary Markov chain,
defined on a probability space $(\Omega,\mathcal{F},P)$ with values in a
Polish space $(S,\mathcal{A})$. Denote by $\mathcal{F}_{n}=\sigma(\xi
_{k},k\leq n)$ and by $\mathcal{F}^{n}=\sigma(\xi_{k},k\geq n)$. The marginal
distribution on $\mathcal{A}$ is denoted by $\pi(A)=\mathbb{P}(\xi_{0}\in
A)$.{ }We shall construct the Markov chain in a canonical way on $S^{Z}$ from
a kernel $Q(x,A)$, and we assume that an invariant distribution $\pi$ exists.{
}

{Next, let }$L${$_{0}^{2}(\pi)$ be the set of measurable functions on $S$ such
that $\int f^{2}d\pi<\infty$ and $\int fd{\pi}=0.$ For a function} ${f}\in
L${$_{0}^{2}(\pi)$ let }%
\begin{equation}
{X_{i}=f(\xi_{i}),\ S_{n}=\sum\nolimits_{i=1}^{n}X_{i}.} \label{def X,S}%
\end{equation}
Denote the regular conditional probability on $\mathcal{F}$, with respect to
$\mathcal{F}_{0}$ by
\[
P^{0}(\cdot)(\omega)=P(\cdot|\mathcal{F}_{0})(\omega),
\]
and the conditional expectation, $E^{0}(X)=E(X|\mathcal{F}_{0}).$ By the
Markov property, if $A\in\mathcal{F}^{0}=\sigma(\xi_{i},i\geq0),$ we have
$P^{0}(A)=P(A|\mathcal{\xi}_{0})$, and for $X$ measurable with respect to
$\mathcal{F}^{0},$ $E^{0}(X)=E(X|\mathcal{\xi}_{0}).$ We are studying the
quenched central limit theorem for Markov chains, which can be stated in two
equivalent ways: For $P$-almost all $\omega\in\Omega$
\begin{equation}
P^{0}(S_{n}/\sqrt{n}\leq t)(\omega)\rightarrow P(N(0,\sigma^{2})\leq t)\text{
for any }t, \label{quenched 1}%
\end{equation}
where $N(0,\sigma^{2})$ is a normal random variable with mean $0$ and variance
$\sigma^{2}$.

Another formulation is known under the name of the CLT started at a point. Let
$P^{x}$ be the probability associated to the Markov chain started from $x\in
S$ and $E^{x}$ be the corresponding conditional expectation. Then, for $\pi
$-almost every $x\in S$,
\begin{equation}
P^{x}(S_{n}/\sqrt{n}\leq t)\rightarrow P(N(0,\sigma^{2})\leq t)\text{ for any
}t. \label{quenched 2}%
\end{equation}
Clearly the quenched CLT implies that for any $t$
\begin{equation}
P(S_{n}/\sqrt{n}\leq t)\rightarrow P(N(0,\sigma^{2})\leq t), \label{annealed}%
\end{equation}
where $N(0,\sigma^{2})$ is a normal random variable with mean $0$ and variance
$\sigma^{2}$. This is called annealed CLT. On the other hand, there are
numerous examples of processes satisfying the annealed CLT but failing to
satisfy the quenched CLT. Some examples of this kind have been constructed by
Voln\'{y} and Woodroofe (\cite{VW}, \cite{VW3}). Therefore, some additional
conditions are needed in order for the central limit theorem to hold in the
quenched form.

The limit theorems started at a point are often encountered in evolutions in
random media and they are of considerable importance in statistical mechanics.
They are also useful for analyzing Markov chain Monte Carlo algorithms. Due to
its importance, the problem was intensively studied in the literature. Two of
the most influential papers are due to Derriennic and Lin (\cite{DL1},
\cite{DL2}), which opened the way for many further results we shall mention
throughout the paper. For a survey on quenched invariance principles under
projective conditions we direct to \cite{Pelsurvey}.

The difficulty of obtaining quenched limit theorems consists in the fact that
a Markov chain started at a point is no longer stationary. This is the reason
this problem is very difficult to solve and there are still many open problems
and long standing conjectures to be settled. Since stationary martingales
satisfy the quenched CLT, the best technique to solve such a problem is to
obtain a martingale approximation with a suitable rest. This technique was
successfully used to get quenched CLT's for various classes of random
variables in numerous papers, \cite{DL1}, \cite{DL2}, \cite{WU}, \cite{CP},
\cite{VW}, \cite{VW2}, \cite{MPP}, \cite{CM}, \cite{DMP}, \cite{BPP}, among
others. The novelty here is that we use a martingale construction and
approximation based on a new idea of conditioning with respect to both the
past and the future of the Markov chain. This idea was introduced in
\cite{Pel20a}, and \cite{Pel20b}. In the annealed setting, if a stationary and
ergodic Markov chain satisfies $E(S_{n}^{2})/n$ is convergent, then the
CLT\ holds (pending only a random centering) (see \cite{Pel20a}). By using a
similar martingale construction we obtain in this paper a new almost sure
martingale approximation under $P^{x},$ for $\pi-$almost all starting points.
This approximation will lead to the quenched CLT under the main condition that
$E^{x}(S_{n}^{2})/n$ is convergent $\pi-$almost surely. As application, we
point out a new class of Markov chains satisfying the quenched CLT, defined by
using projective conditions. In defining this class no assumption of
irreducibility nor of aperiodicity is imposed. Under the additional
assumptions that the Markov chain is irreducible, aperiodic and positively
recurrent, Chen (Proposition 3.1., \cite{Chen}) showed that if the CLT holds
for the stationary Markov chain then the quenched CLT holds.

Here are some notations we shall use throughout the paper. {We denote by
}${{||X||}}$ the norm in $L${$^{2}$}$(\Omega,\mathcal{F},P)$.{ }Unless
otherwise specified, we shall assume the total ergodicity of the shift $T$ of
the sequence $(\xi_{n})_{n\in\mathbb{Z}}$ with respect to $P,$ i.e. $T^{m}$ is
ergodic for every $m\geq1$. For the definition of the ergodicity of the shift
we direct the reader to the subsection "A return to Ergodic Theory" in
Billingsley \cite{Bil} p. 494. Let us consider the operator $Q$ induced by the
kernel $Q(x,A)$ on bounded measurable functions on $(S,\mathcal{A})$ defined
by $Qf(x)=\int\nolimits_{S}f(y)Q(x,dy)$. By using Corollary 5 p. 97 in
Rosenblatt \cite{Ros1}, the shift of $(\xi_{n})_{n\in\mathbb{Z}}$ is totally
ergodic with respect to $P$\ if and only if the powers $Q^{m}\ $are ergodic
with respect to $\pi$ for all natural $m$ (i.e. $Q^{m}f=f$ for $f$ bounded on
$(S,\mathcal{A})$ implies $f$ is constant $\pi-$a.e.). For more information on
total ergodicity, we refer to the survey paper by Quas \cite{Quas}.

Throughout the paper $\Rightarrow$ denotes the convergence in distribution. By
the notation a.s. we understand $P$-almost surely. We shall also use the
notation $K$ for the conditional expectation operator on $L_{1}(P)$, namely
\[
K(X)=E(X\circ T^{-1}|\mathcal{\xi}_{0}),\text{ \ }K^{n}(X)=K(K^{n-1}%
(X))=E(E(X\circ T^{-n})|\mathcal{\xi}_{0}).
\]

The problem we address in Theorem \ref{Th Main} is to provide necessary and
sufficient conditions for a quenched CLT for a class of Markov chains.

\begin{theorem}
\label{Th Main} Assume $(X_{n})$ and $(S_{n})$ are defined by (\ref{def X,S}),%
\begin{equation}
\lim\sup_{n\rightarrow\infty}\frac{E(S_{n}^{2})}{n}<\mathbb{\infty}
\label{varsup1}%
\end{equation}
and
\begin{equation}
\lim_{n\rightarrow\infty}\frac{1}{n}||E(S_{n}|\xi_{0},\xi_{n})||^{2}=0.\text{
} \label{neglipf}%
\end{equation}
Then there is $\sigma\geq0$ such that%
\begin{equation}
\frac{S_{n}}{\sqrt{n}}\Rightarrow N(0,\mathbb{\sigma}^{2})\text{ and }%
\frac{E(S_{n}^{2})}{n}\rightarrow\sigma^{2}. \label{CLT+moments}%
\end{equation}
Furthermore, the following are equivalent:%
\[
(a)\text{
\ \ \ \ \ \ \ \ \ \ \ \ \ \ \ \ \ \ \ \ \ \ \ \ \ \ \ \ \ \ \ \ \ \ \ \ \ }%
\lim\sup_{n\rightarrow\infty}\frac{E^{0}\left(  S_{n}^{2}\right)  }{n}%
<\sigma^{2}\text{ }\ \text{a.s.
\ \ \ \ \ \ \ \ \ \ \ \ \ \ \ \ \ \ \ \ \ \ \ \ \ \ \ \ }%
\]%
\[
(b)\text{
\ \ \ \ \ \ \ \ \ \ \ \ \ \ \ \ \ \ \ \ \ \ \ \ \ \ \ \ \ \ \ \ \ \ \ \ \ \ }%
\frac{E^{0}(S_{n}^{2})}{n}\text{ converges }\ \text{a.s.
\ \ \ \ \ \ \ \ \ \ \ \ \ \ \ \ \ \ \ \ \ \ \ \ \ \ \ \ \ \ }%
\]
$(c)$ The quenched CLT in (\ref{quenched 1})\ holds and $S_{n}^{2}/n$ is
uniformly integrable under $P^{0}(\omega)$ for almost all $\omega.$
\end{theorem}

\begin{remark}
Note that in condition (b) we do not have to specify the almost sure limit of
$E^{0}(S_{n}^{2})/n.$ However, under our conditions it will always be
$\sigma^{2}$. In the sequel, when we say that the quenched limit theorem holds
we understand that (c) of Theorem \ref{Th Main} holds.
\end{remark}

Relevant for the next Corollary is the notion of two-sided tail sigma field.
We define the two-sided tail sigma field by
\[
\mathcal{T}_{d}=\cap_{n\geq1}(\mathcal{F}_{-n}\vee\mathcal{F}^{n}).
\]
We say that $\mathcal{T}_{d}$ is trivial if for any $A\in\mathcal{T}_{d}$ we
have $P(A)=0$ or $1$.

\begin{corollary}
\label{corr tail}Assume $(X_{n})$ and $(S_{n})$ are defined by (\ref{def X,S}%
), $\mathcal{T}_{d}$ is trivial and $S_{n}^{2}/n$ is uniformly integrable.
Then the Markov chain is totally ergodic, (\ref{CLT+moments}) holds and in
addition (a), (b), and (c) of Theorem \ref{Th Main} are equivalent.
\end{corollary}

In the next section we shall point out a sufficient condition for the quenched
CLT by using projective criteria.

\section{A sufficient condition for the quenched CLT}

In this section we give a new sufficient condition for the quenched CLT based
on the proof of Theorem \ref{Th Main}. This condition arises in a computation
of $E^{0}\left(  S_{n}^{2}\right)  $ by dyadic expansion.

We recall that the sequences $(X_{n})$ and $(S_{n})$ are defined by
(\ref{def X,S}).

As shown in Theorem 2.7 in Cuny and Merlev\`{e}de\ \cite{CM}, it is known that
the quenched CLT\ holds under a condition introduced by Maxwell and Woodroofe
\cite{MW}, namely%
\begin{equation}
\sum_{n\geq1}\frac{||E(S_{n}|\xi_{0})||}{n^{3/2}}<\infty. \label{cond MW}%
\end{equation}

There are examples of Markov chains pointing out that, in general, condition
(\ref{cond MW}) is as sharp as possible in some sense. Peligrad and Utev
\cite{PU} constructed an example showing that for any sequence of positive
constants $(a_{n}),$ $a_{n}\rightarrow0,$ there exists a stationary Markov
chain such that
\[
\sum_{n\geq1}a_{n}\frac{||E(S_{n}|\xi_{0})||}{n^{3/2}}<\infty
\]
but $S_{n}/\sqrt{n}$ is not stochastically bounded. This example and other
counterexamples provided by Voln\'{y} \cite{V}, Dedecker \cite{Ded} and Cuny
and Lin \cite{CL}, show that, in general, condition
\begin{equation}
\sum_{n\geq1}\frac{||E(S_{n}|\xi_{0})||^{2}}{n^{2}}<\infty\label{conjrev}%
\end{equation}
does not assure that $(S_{n}/\sqrt{n})$ is stochastically bounded. However,
Corollary 3.5 in \cite{Pel20b} contains a CLT under a reinforced form of
(\ref{conjrev}).\ \ We provide next a quenched form of that result.

\begin{theorem}
\label{ThPeligrad}The quenched CLT holds under the condition
\begin{equation}
\sum_{n\geq1}\frac{||E(S_{n}|\xi_{0},\xi_{n})||^{2}}{n^{2}}<\infty.
\label{condpf}%
\end{equation}

\end{theorem}

As a corollary to Theorem \ref{ThPeligrad}, by Lemma 14 in \cite{Pel20b} we
have the following sufficient condition for (\ref{condpf}) in terms of
individual summands:

\begin{corollary}
The quenched CLT holds under the condition%
\begin{equation}
\sum\nolimits_{k\geq1}||E(X_{0}|\xi_{-k},\xi_{k})||^{2}<\infty.
\label{mixingale}%
\end{equation}

\end{corollary}

We end this section by mentioning a conjecture due to Kipnis and Varadhan
\cite{KV}, which is unsolved. The conjecture asks if the quenched CLT and its
functional form hold for stationary reversible and ergodic Markov chains
($Q=Q^{\ast}$ with $Q^{\ast}$ the adjoint of $Q$) satisfying (\ref{conjrev}).
For reversible Markov chains (\ref{conjrev}) is an equivalent formulation of
$E(S_{n}^{2})/n$ converges. This problem was investigated in several papers,
\cite{DL1}, \cite{CP} where the quenched CLT for reversible Markov chains was
obtained under various reinforcements of (\ref{conjrev}).

\section{Proofs}

The starting point of the proofs is a new annealed CLT for Markov chains (see
Theorem 1 in \cite{Pel20a}):

\begin{theorem}
\label{Th random center}Let $(X_{n})_{n\in Z}$ and $(S_{n})_{n\geq1}$ be as
defined in (\ref{def X,S}), $(\xi_{n})$ is totally ergodic, and assume that
(\ref{varsup1}) holds. Then, the following limit exists
\begin{equation}
\lim_{n\rightarrow\infty}\frac{1}{n}||S_{n}-E(S_{n}|\xi_{0},\xi_{n}%
)||^{2}=\sigma^{2} \label{defc}%
\end{equation}
and%
\[
\frac{S_{n}-E(S_{n}|\xi_{0},\xi_{n})}{\sqrt{n}}\Rightarrow N(0,\sigma
^{2})\text{ as }n\rightarrow\infty.
\]

\end{theorem}

This result has the following consequence: (Corollary 5, \cite{Pel20a}):

\begin{theorem}
\label{Th CLT}Assume that (\ref{varsup1})\ and (\ref{neglipf})\ holds. Then
(\ref{CLT+moments}) holds.
\end{theorem}

The main step of proving Theorem \ref{Th Main} is the following proposition:

\begin{proposition}
\label{pr quenched}If in addition to the conditions of Theorem \ref{Th CLT} we
assume that
\begin{equation}
\lim\sup_{n\rightarrow\infty}\frac{E^{0}(S_{n}^{2})}{n}\leq\mathbb{\sigma}%
^{2}\text{ }\ \text{a.s.} \label{cond bound}%
\end{equation}
then the quenched CLT\ in (\ref{quenched 1}) holds.
\end{proposition}

\bigskip

\textbf{Proof of Proposition \ref{pr quenched}}

\bigskip

The proof of the quenched CLT is also based on the new idea to use a
martingale approximation by conditioning with respect to past and future of
the chain. We shall use the notations $E(X^{2}|\xi_{0},\xi_{n})=||X||_{0,n}%
^{2}$ and $E^{0}(X^{2})=||X||_{0}^{2}.$

We start the proof by a decomposition in blocks of random variables, which is
intended to weaken the dependence. Fix $m$ ($m<n$) a positive integer and make
consecutive blocks of size $m$. Denote by $Y_{k}$ the sum of variables in the
$k$'th block. Let $u=u_{n}(m)=[n/m].$ So, for $k=0,1,...,u-1,$ we have%
\begin{equation}
Y_{k}=Y_{k}(m)=(X_{km+1}+...+X_{(k+1)m}). \label{defY}%
\end{equation}
Also denote%
\[
Y_{u}=Y_{u}(m)=(X_{um+1}+...+X_{n}).
\]
With this notations we write%
\[
\frac{1}{\sqrt{u}}S_{u}(m):=\frac{1}{\sqrt{u}}\sum\nolimits_{k=0}^{u-1}%
\frac{1}{\sqrt{m}}Y_{k}(m)=\frac{1}{\sqrt{um}}S_{mu}.
\]
In the first step of the proof we show that it is enough to prove that
$S_{u}(m)/\sqrt{u}$ satisfies the quenched CLT. Let us show that the last
block $Y_{u}(m)/\sqrt{n}$ has a negligible contribution to the convergence in
distribution. With this aim, by Theorem 3.1 in Billingsley \cite{Bil}, it is
enough to show that
\begin{equation}
E^{0}\left(  \frac{S_{n}-S_{mu}}{\sqrt{n}}\right)  ^{2}=E^{0}\left(
\frac{Y_{u}(m)}{\sqrt{n}}\right)  ^{2}\rightarrow0\text{ a.s. as }%
n\rightarrow\infty. \label{negli3}%
\end{equation}
Note that the definition of $Y_{u}(m)$ and the Cauchy-Schwartz inequality
imply that
\[
E^{0}\left(  \frac{Y_{u}(m)}{\sqrt{n}}\right)  ^{2}\leq m\frac{\max_{1\leq
j\leq n}E^{0}(X_{j}^{2})}{n}.
\]
Now, fix $M>0$ and note that, for each $\varepsilon>0$ and $n>M,$
\begin{align*}
\frac{\max_{1\leq j\leq n}E^{0}(X_{j}^{2})}{n}  &  \leq\varepsilon^{2}%
+\frac{{\sum\nolimits_{j=1}^{n}}E^{0}(X_{j}^{2}{I(|X_{j}|>\varepsilon}\sqrt
{n}))}{n}\\
&  \leq\varepsilon^{2}+\frac{{\sum\nolimits_{j=1}^{n}E^{0}(X_{j}^{2}%
I(|X_{j}|>\varepsilon}\sqrt{M}{)}}{n}.
\end{align*}
So, by Hopf's pointwise ergodic theorem for Dunford--Schwartz operators
(Theorem 7.3 in Krengel \cite{Kr})
\[
\lim\sup_{n\rightarrow\infty}\frac{\max_{1\leq j\leq n}E^{0}(X_{j}^{2})}%
{n}\leq\varepsilon^{2}+E({X_{0}^{2}I(|X_{0}|>\varepsilon}\sqrt{M}{)}\text{
a.s.}%
\]
and so, letting $\varepsilon\rightarrow0$ and $M\rightarrow\infty$ we have
\[
\lim\sup_{n\rightarrow\infty}\frac{\max_{1\leq j\leq n}E^{0}(X_{j}^{2})}%
{n}=0\text{ a.s.}%
\]

By the above arguments, we have proved that (\ref{negli3}) holds for any $m$,
and therefore $S_{n}/\sqrt{n}$ has the same limiting distribution as
$S_{um}/\sqrt{n}$ under $P^{0}(\omega)$ for almost all $\omega$. Since
$um/n=[n/m](m/n)\rightarrow1$ as $n\rightarrow\infty,$ by Slutsky's theorem,
$S_{um}/\sqrt{n}$ has the same limiting distribution as $S_{u}(m)/\sqrt{u}.$
Furthermore, from (\ref{cond bound}) and (\ref{negli3}) we easily derive that%
\begin{equation}
\lim\sup_{u\rightarrow\infty}\frac{1}{u}||S_{u}(m)||_{0\ }^{2}\leq\sigma
^{2}\text{ a.s.} \label{cond bound 2}%
\end{equation}
In the second step of the proof we construct the approximating martingale and
mention its limiting properties.

For $k=0,1,...,u-1,$ let us consider the random variables\textbf{ }%
\[
D_{k}=D_{k}(m)=\frac{1}{\sqrt{m}}(Y_{k}-E(Y_{k}|\xi_{km},\xi_{(k+1)m})).
\]
By the Markov property, conditioning by $\sigma(\xi_{km},\xi_{(k+1)m})$ is
equivalent to conditioning by $\mathcal{F}_{km}\vee\mathcal{F}^{(k+1)m}.$ Note
that $D_{k}$ is adapted to $\mathcal{F}_{(k+1)m}=\mathcal{G}_{k}.$ Also note
that $\mathcal{G}_{0}=\sigma(\xi_{i},i\leq0).$ Then we have $E(D_{1}%
|\mathcal{G}_{0})=0$ a.s. Since we assumed that the shift $T$ of the sequence
$(\xi_{n})_{n\in\mathbb{Z}}$ is totally ergodic, we deduce that for every $m$
fixed, we have a stationary and ergodic sequence of square integrable
martingale differences $(D_{k},\mathcal{G}_{k})_{k\geq0}$.

Therefore, by the classical quenched central limit theorem for ergodic
martingales, (see page 520 in Derriennic and Lin \cite{DL1}) for every $m,$ a
fixed positive integer, we have for almost all $\omega\in\Omega$,
\[
\frac{1}{\sqrt{u}}M_{u}(m):=\frac{1}{\sqrt{u}}\sum\nolimits_{k=0}^{u-1}%
D_{k}(m)\Rightarrow N_{m}\text{ as }u\rightarrow\infty,\text{ under }%
P^{0}(\omega)\text{,}%
\]
where $N_{m}$ is a normally distributed random variable with mean $0$ and
variance $E(D_{0}^{2})=m^{-1}||S_{m}-E(S_{m}|\xi_{0},\xi_{m})||^{2}$.

Since by (\ref{neglipf}) and (\ref{defc}),
\begin{equation}
m^{-1}||S_{m}-E(S_{m}|\xi_{0},\xi_{m})||^{2}=m^{-1}(||S_{m}||^{2}%
-||E(S_{m}|\xi_{0},\xi_{m})||^{2})\rightarrow\sigma^{2}, \label{var mart}%
\end{equation}
it follows that $N_{m}\Rightarrow N(0,\sigma^{2})$. So, for almost all
$\omega\in\Omega,$%
\begin{equation}
\frac{1}{\sqrt{u}}M_{u}(m)\Rightarrow N_{m}\Rightarrow N(0,\sigma^{2})\text{
under }P^{0}(\omega)\text{.} \label{star}%
\end{equation}
In the last step of the proof we shall approximate $S_{u}(m)$ by $M_{u}(m)$ in
a suitable way, which will allow us to get the quenched limiting distribution
$N(0,\sigma^{2})$ also for $S_{u}(m)/\sqrt{u},$ completing the proof of this
theorem. By using Theorem 3.2 in Billingsley \cite{bil2} and taking into
account (\ref{star}),\ in order to establish the quenched CLT from Proposition
\ref{pr quenched}, we have only to show that
\begin{equation}
\lim\inf_{m\rightarrow\infty}\lim\sup_{u\rightarrow\infty}E^{0}(\frac{1}%
{\sqrt{u}}S_{u}(m)-\frac{1}{\sqrt{u}}M_{u}(m))^{2}=0\text{ a.s.}
\label{negli1}%
\end{equation}
Denote by
\begin{equation}
Z_{k}=m^{-1/2}E(Y_{k}|\xi_{km},\xi_{(k+1)m})\text{ and }R_{u}(m)=\sum
\nolimits_{k=0}^{u-1}Z_{k}. \label{defRn}%
\end{equation}
With this notation we have:
\begin{equation}
S_{u}(m)=M_{u}(m)+R_{u}(m). \label{decom}%
\end{equation}
Let us show that $M_{u}(m)$ and $R_{u}(m)$ are orthogonal given $\mathcal{F}%
_{0}\vee\mathcal{F}^{n}$. We show this property by analyzing the conditional
expected value of all the terms of the product $M_{u}(m)R_{u}(m)$. For $m\leq
n,$ and $X\in\sigma(\xi_{j},$ $m\leq j\leq n)$ it is convenient to denote
$E^{m,n}(X)=E(X|\mathcal{F}_{m}\vee\mathcal{F}^{n})=E(X|\xi_{m}\vee\xi_{n})$.
Note that if $j<k,$ since $\mathcal{F}_{(j+1)m}\subset\mathcal{F}_{km}$, and
taking into account the Markov chain properties, we have that
\begin{align*}
&  E^{0,n}[(Y_{k}-E(Y_{k}|\xi_{km},\xi_{(k+1)m}))E(Y_{j}|\xi_{jm},\xi
_{(j+1)m})]\\
&  =E^{0,n}[E^{(j+1)m,n}(Y_{k}-E(Y_{k}|\xi_{km},\xi_{(k+1)m}))E(Y_{j}|\xi
_{jm},\xi_{(j+1)m})]\\
&  =E^{0,n}[E^{(j+1)m,n}(Y_{k}-E(Y_{k}|\mathcal{F}_{km}\vee\mathcal{F}%
^{(k+1)m}))E(Y_{j}|\xi_{jm},\xi_{(j+1)m})]=0\text{ a.s.}%
\end{align*}
On the other hand, if $j>k,$ since $\mathcal{F}^{jm}\subset\mathcal{F}%
^{(k+1)m}$ then
\begin{align*}
&  E^{0,n}[(Y_{k}-E(Y_{k}|\xi_{km},\xi_{(k+1)m}))E(Y_{j}|\xi_{jm},\xi
_{(j+1)m})]\\
&  =E^{0,n}[E^{0,jm}(Y_{k}-E(Y_{k}|\xi_{km},\xi_{(k+1)m}))E(Y_{j}|\xi_{jm}%
,\xi_{(j+1)m})]\\
&  =E^{0,n}[E^{0,jm}(Y_{k}-E(Y_{k}|\mathcal{F}_{km}\vee\mathcal{F}%
^{(k+1)m}))E(Y_{j}|\xi_{jm},\xi_{(j+1)m})]=0\text{ a.s.}%
\end{align*}
For $j=k,$ by conditioning with respect to $\sigma(\xi_{km},\xi_{(k+1)m}),$ we
note that
\[
E^{0,n}[(Y_{k}-E(Y_{k}|\xi_{km},\xi_{(k+1)m}))E(Y_{k}|\xi_{km},\xi
_{(k+1)m})]=0\text{ a.s.}%
\]
Therefore $M_{u}(m)$ and $R_{u}(m)$ are indeed orthogonal under $E^{0,n}$
almost surely. By using now the decomposition (\ref{decom}), and the fact that
$M_{u}(m)$ and $R_{u}(m)$ are orthogonal a.s. under $E^{0,n}$, we obtain the
identity
\begin{equation}
\frac{1}{u}||S_{u}(m))||_{0,n}^{2}=\frac{1}{u}||M_{n}(m)||_{0,n}^{2}+\frac
{1}{u}||R_{u}(m)||_{0,n}^{2}\text{ a.s.} \label{ID}%
\end{equation}
By conditioning with respect to $\sigma(\xi_{0})$ in (\ref{ID}), and taking
into account the properties of conditional expectation, we also have%
\[
\frac{1}{u}||S_{u}(m)||_{0}^{2}=\frac{1}{u}||M_{n}(m)||_{0}^{2}+\frac{1}%
{u}||R_{u}(m)||_{0}^{2}\text{.}%
\]
By the definition of $M_{u}(m)$,
\[
\frac{1}{u}||M_{n}(m)||_{0}^{2}=\frac{1}{u}\sum\nolimits_{k=0}^{u-1}E^{0}%
D_{k}^{2}(m)=\frac{1}{u}\sum\nolimits_{k=0}^{u-1}K^{k}\left(  D_{0}%
^{2}(m)\right)  \text{.}%
\]
Now, by using the fact that $(\xi_{n})$ is totally ergodic along with Hopf's
pointwise ergodic theorem for Dunford--Schwartz operators,%

\[
\lim_{u\rightarrow\infty}\frac{1}{u}||M_{u}(m)||_{0}^{2}=\frac{1}{m}%
||S_{m}-E(S_{m}|\xi_{0},\xi_{m})||^{2}\text{ a.s.}%
\]
So, by (\ref{var mart})%
\[
\lim_{m\rightarrow\infty}\lim_{u\rightarrow\infty}\frac{1}{u}||M_{u}%
(m)||_{0}^{2}=\sigma^{2}.
\]
By passing now to the limit in (\ref{ID}) and using (\ref{cond bound 2}) we
obtain%
\[
\sigma^{2}\geq\lim\sup_{u\rightarrow\infty}\frac{1}{u}||S_{u}(m)||_{0\ }%
^{2}\geq\sigma^{2}+\lim\sup_{m\rightarrow\infty}\lim\sup_{u\rightarrow\infty
}\frac{1}{u}||R_{u}(m)||_{0}^{2}\text{ \ a.s.}%
\]
Therefore,%
\[
\lim_{m\rightarrow\infty}\lim\sup_{u\rightarrow\infty}\frac{1}{u}%
||R_{u}(m)||_{0}^{2}=0\text{ \ a.s.,}%
\]
which implies (\ref{negli1}),\ and also the result follows. $\square$

\bigskip

In the next lemma we mention a property of the limit of $E^{0}(S_{n}^{2})/n.$
The idea of proof is borrowed from Dedecker and Merlev\`{e}de \cite{DM},
Subsection (3.2), where it was used in another context.

\begin{lemma}
\label{Invariance}Assume that%
\[
\frac{1}{n}E^{0}(S_{n}^{2})\rightarrow\eta\text{ in }L_{1}.
\]
Then $\eta$ is measurable with respect to the invariant sigma field.
\end{lemma}

Proof. Recall the definition of shift $T$. Below, we denote by $TX=X\circ
T^{-1}.$

Clearly, $\eta$ is $\mathcal{F}_{0}$ measurable. Then
\begin{equation}
E\left\vert E^{0}\left(  \frac{1}{n}S_{n}^{2}-\eta\right)  \right\vert
\rightarrow0\text{ as }n\rightarrow\infty. \label{rel1}%
\end{equation}
Therefore, with the notation $E^{1}(\cdot)=E(\cdot|\mathcal{F}_{1}),$%

\[
E\left\vert E^{1}\left(  \frac{1}{n}TS_{n}^{2}-T\eta\right)  \right\vert
\rightarrow0\text{ as }n\rightarrow\infty.
\]
Since $\mathcal{F}_{0}\subset\mathcal{F}_{1},$ by the properties of
conditional expectation, this implies%
\[
E\left\vert E^{0}\left(  \frac{1}{n}TS_{n}^{2}-T\eta\right)  \right\vert
\rightarrow0\text{ as }n\rightarrow\infty.
\]
But, since the condition of this lemma implies that $E(S_{n}^{2})/n$ is
bounded,
\begin{align*}
\frac{1}{n}E|S_{n}^{2}-TS_{n}^{2}|  &  \leq\frac{1}{n}E|(S_{n}^{2}%
-(S_{n}-X_{1}+X_{n+1})^{2}|\\
&  \leq\frac{1}{n}E|(X_{1}-X_{n+1})(2S_{n}-X_{1}+X_{n+1})|\\
&  \leq\frac{4}{n}||X_{0}||\cdot\left(  ||S_{n}||+||X_{0}||\right)
\rightarrow0\text{ as }n\rightarrow\infty.
\end{align*}
So, by combining the last two limits, we also have
\[
E\left\vert E^{0}\left(  \frac{1}{n}S_{n}^{2}-T\eta\right)  \right\vert
\rightarrow0\text{ as }n\rightarrow\infty.
\]
By combining this limit with (\ref{rel1}) we obtain%
\[
E|\left(  \eta-E^{0}\left(  T\eta\right)  \right)  |=0,
\]
implying that
\[
\eta=E^{0}\left(  T\eta\right)  \text{ a.s.}%
\]
It remains to apply Lemma 3 from Dedecker and Merlev\`{e}de \cite{DM}, giving
that $\eta=T\eta$ a.s. $\square$

\bigskip

\textbf{Proof of Theorem \ref{Th Main}}

\bigskip

The first part in Theorem \ref{Th Main} is given in Theorem \ref{Th CLT}, so
we have to prove only the second part of this theorem.

We argue first that (a) implies (c).

Since we assume (a) the quenched CLT holds by Proposition \ref{pr quenched}.
Note that, by Theorem 25.11 in Billingsley \cite{bil2}, the quenched CLT
implies that
\[
\sigma^{2}\leq\lim\inf_{n\rightarrow\infty}E^{0}(S_{n}^{2})/n\text{ a.s.,}%
\]
which combined with (a) gives $E^{0}(S_{n}^{2})/n\rightarrow\sigma^{2}$ a.s.
Now, because we have the quenched CLT and $E^{0}(S_{n}^{2})/n\rightarrow
\sigma^{2}$ a.s., by Theorem 3.6 in Billingsley \cite{bil2}, we have the
uniform integrability of $(S_{n}^{2}/n)_{n}$ under $P^{0}(\omega)$ for almost
all $\omega.$

Clearly (c) implies (b) by the convergence of the moments in the CLT in
Theorem 3.5 in Billingsley \cite{bil2}. Actually (c) implies $E^{0}(S_{n}%
^{2})/n\rightarrow\sigma^{2}$ a.s.

It remains to show that (b) implies (a).

We start from (b), which is: for some random variable $\eta$, $E^{0}(S_{n}%
^{2})/n\rightarrow\eta$ a.s. By Theorem \ref{Th CLT} the annealed
CLT\ together with the convergence of the second moments hold. Furthermore, by
Theorem 3.6 Billingsley \cite{bil2} we have that $S_{n}^{2}/n$ is uniformly
integrable. This implies that $E^{0}(S_{n}^{2})/n$ is also uniformly
integrable, which, together with (b), implies the convergence $E^{0}(S_{n}%
^{2})/n\rightarrow\eta$ in $L^{1}(\Omega,\mathcal{F},P).$ By Lemma
\ref{Invariance}, the limit of $E^{0}(S_{n}^{2})/n$ is measurable with respect
to the trivial invariant sigma field; therefore it is constant. Because we
assumed that $E(S_{n}^{2})/n\rightarrow\sigma^{2}$ it follows that
$\eta=\sigma^{2}.$ $\square$

\bigskip

\textbf{Proof of Corollary \ref{corr tail}}

\bigskip

The fact that $(\xi_{k})_{k\in Z}$ is totally ergodic follows by Proposition
2.12 in the Vol. 1 of Bradley (2007). Then, since we assume that $\left(
S_{n}^{2}/n\right)  $ is uniformly integrable, by Lemma 4 in \cite{Pel22} we
deduce that (\ref{neglipf}) holds. The result follows by Theorem
\ref{Th Main}. $\square$

\bigskip

We move now to prove Theorem \ref{ThPeligrad}. Relevant for the proof is the
weak $L_{p}$ space, define by
\[
L^{p,w}=\{f\text{ measurable, }\sup_{\lambda>0}\lambda^{p}P(|f|\geq
\lambda)<\infty\}.
\]
Denote the norm in $L^{p,w}$ by $||\cdot||_{p,w}.$ Below we also use the
notation $\bar{S}_{k}=S_{2k}-S_{k}.$

The main step for proving Theorem \ref{ThPeligrad} is the following upper
bound concerning $E^{0}\left(  S_{n}^{2}\right)  /n.$

\begin{lemma}
\label{lmbound var}For any stationary and ergodic sequence $(\eta_{n})$, not
necessarily Markov, define $(V_{n})$ by $V_{n}=g(\eta_{n})$ and $S_{n}%
=\sum\nolimits_{k=1}^{n}V_{k}.$ Assume $V_{0}$ is in $L_{2}$ and is centered
at expectation. Let $\mathcal{K}_{n}=\sigma(\eta_{j},j\leq n)$ and keep the
notation $E^{n}(X)=E(X|\mathcal{K}_{n}).$ Then we have the following bound%
\[
||\sup_{n}\frac{1}{n}E^{0}\left(  S_{n}^{2}\right)  ||_{1,w}\leq6E(V_{0}%
^{2})+12\sum\nolimits_{k\geq0}\frac{1}{2^{k}}E|E^{0}\left(  S_{2^{k}}\bar
{S}_{2^{k}}\right)  |\text{ .}%
\]

\end{lemma}

Proof. The proof follows the traditional technique of dyadic recurrence,
initiated by Ibragimov \cite{Ib} and further developed in \cite{P}, \cite{Br},
\cite{PU}, \cite{CM}, among many others.

Let $2^{r-1}\leq n<2^{r}$ and write its binary expansion:
\[
n=\sum_{k=0}^{r-1}2^{k}a_{k}\text{ where }a_{r-1}=1\text{ and }a_{k}%
\in\{0,1\}\text{ for }k=0,\dots,r-2\,.
\]
Notice that
\[
S_{n}=\sum_{i=0}^{r-1}a_{i}T_{2^{i}}\text{ where $T_{2^{i}}=\sum_{i=n_{i-1}%
+1}^{n_{i}}X_{i}$, $\ n_{i}=\sum_{j=0}^{i}a_{j}2^{j}$ and $n_{-1}=0$%
}\,\text{.}%
\]
By the triangle inequality, (recall that $S_{0}=0$)%
\[
(E^{0}(S_{n}^{2}))^{1/2}=||S_{n}||_{0}=||\sum_{i=0}^{r-1}a_{i}E^{0}(S_{n_{i}%
}-S_{n_{i-1}})||_{0}\leq\sum_{i=0}^{r-1}||S_{n_{i}}-S_{n_{i-1}}||_{0}.
\]
Also, by stationarity and because $n_{i}-n_{i-1}$ is either $0$ or $2^{i},$ we
obtain%
\begin{align*}
E^{0}(S_{n_{i}}-S_{n_{i-1}})^{2}  &  =E^{0}(E((S_{n_{i}}-S_{n_{i-1}}%
)^{2}|\mathcal{K}_{n_{i-1}}))=K^{n_{i-1}}(E^{0}(S_{n_{i}-n_{i-1}}^{2}))\\
&  \leq K^{n_{i-1}}\left(  E^{0}(S_{2^{i}}^{2})\right)  .
\end{align*}
It follows that%
\begin{equation}
\frac{1}{n}E^{0}(S_{n}^{2})\leq\frac{1}{n}\left(  \sum_{i=0}^{r-1}\left[
K^{n_{i-1}}\left(  E^{0}(S_{2^{i}}^{2})\right)  \right]  ^{1/2}\right)
^{2}\leq6\sup_{i}\frac{1}{2^{i}}\left[  K^{n_{i-1}}\left(  E^{0}(S_{2^{i}}%
^{2})\right)  \right]  . \label{var0}%
\end{equation}
We fix $\ i\geq1\ $ and evaluate the term in the right hand side of
(\ref{var0}).

For each $k$ and $j$, denote $S_{1,2^{k}}=S_{2^{k}},$ $S_{j,2^{k}}=S_{j2^{k}%
}-S_{(j-1)2^{k}}$. Clearly,%
\begin{align*}
S_{2^{i}}^{2}  &  =S_{1,2^{i-1}}^{2}+S_{2,2^{i-1}}^{2}+2S_{1,2^{i-1}%
}S_{2,2^{i-1}}=S_{1,2^{i-2}}^{2}+S_{2,2^{i-2}}^{2}+S_{3,2^{i-2}}%
^{2}+S_{4,2^{i-2}}^{2}\\
&  +2\left(  S_{1,2^{i-2}}S_{2,2^{i-2}}+S_{3,2^{i-2}}S_{4,2^{i-2}%
}+S_{1,2^{i-1}}S_{2,2^{i-1}}\right)  .
\end{align*}
We continue the recurrence and get the representation:%
\[
S_{2^{i}}^{2}=\sum\nolimits_{j=1}^{2^{i}}V_{j}^{2}+2\sum\nolimits_{k=0}%
^{i-1}\sum\nolimits_{j=1}^{2^{i-k-1}}S_{2j-1,2^{k}}S_{2j,2^{k}}.
\]
Denoting by
\[
g_{k}=E^{0}\left(  S_{2^{k}}\bar{S}_{2^{k}}\right)  ,
\]
note that, by using the definition of the conditional expectation $K,$
\[
E^{0}\left(  S_{2j-1,2^{k}}S_{2j,2^{k}}\right)  =E^{0}(E(S_{2j-1,2^{k}%
}S_{2j,2^{k}}|\mathcal{K}_{(2j-2)2^{k}}))=K^{(j-1)2^{k+1}}(g_{k}).
\]
By the above considerations,
\[
\frac{1}{2^{i}}E^{0}(S_{2^{i}}^{2})=\frac{1}{2^{i}}\sum\nolimits_{j=1}^{2^{i}%
}K^{j}(V_{0}^{2})+2\sum\nolimits_{k=0}^{i-1}\frac{1}{2^{i-k}}\left(
\sum\nolimits_{j=1}^{2^{i-k}-1}K^{(j-1)2^{k+1}}\right)  \left(  \frac{1}%
{2^{k}}g_{k}\right)  .
\]
So,%
\begin{align*}
\frac{1}{2^{i}}K^{n_{i-1}}\left(  E^{0}(S_{2^{i}}^{2})\right)   &  =\frac
{1}{2^{i}}\sum\nolimits_{j=1}^{2^{i}}K^{j}\left(  K^{n_{i-1}}(V_{0}%
^{2}\right)  )\\
&  +2\sum\nolimits_{k=0}^{i-1}\frac{1}{2^{i-k}}\left(  \sum\nolimits_{j=1}%
^{2^{i-k}-1}K^{(j-1)2^{k+1}}\right)  \left(  \frac{1}{2^{k}}K^{n_{i-1}}%
(g_{k})\right)  .
\end{align*}
So, with the notation%
\[
\sup_{n}\frac{1}{n}\left(  \sum\nolimits_{j=0}^{n-1}K^{j2^{k+1}}%
(\cdot)\right)  =\mathcal{M}_{k}(\cdot),
\]
we obtain%
\[
\frac{1}{2^{i}}K^{n_{i-1}}\left(  E^{0}S_{2^{i}}^{2}\right)  \leq\sup_{n}%
\frac{1}{n}\sum\nolimits_{j=1}^{n}K^{j}\left(  K^{n_{i-1}}V_{0}^{2}\right)
+2\sum\nolimits_{k=0}^{i-1}\frac{1}{2^{k}}\mathcal{M}_{k}\left(  \left\vert
K^{n_{i-1}}(g_{k})\right\vert \right)  .
\]
By using now Hopf's ergodic theorem (see, e.g., Krengel \cite{Kr}, Lemma 6.1,
page 51, and Corollary 3.8, page 131),%
\[
||\mathcal{M}_{k}\left(  \left\vert K^{n_{i-1}}\left(  \frac{g_{k}}{2^{k}%
}\right)  \right\vert \right)  ||_{1,w}\leq\frac{1}{2^{k}}\left\vert
\left\vert K^{n_{i-1}}g_{k}\right\vert \right\vert _{1}\leq\frac{1}{2^{k}%
}||g_{k}||_{1}\text{.}%
\]
and also%
\[
||\sup_{n}\frac{1}{n}\sum\nolimits_{j=1}^{n}K^{j}\left(  K^{n_{i-1}}(V_{0}%
^{2}\right)  )||_{1,w}\leq E(K^{n_{i-1}}V_{0}^{2})=E(V_{0}^{2}).
\]
Therefore,%
\[
||\sup_{i}\frac{1}{2^{i}}K^{n_{i-1}}\left(  E^{0}(S_{2^{i}}^{2})\right)
||_{1,w}\leq E(V_{0}^{2})+2\sum\nolimits_{k\geq0}\frac{1}{2^{k}}E|E^{0}\left(
S_{2^{k}}\bar{S}_{2^{k}}\right)  |.
\]
To obtain the conclusion of this lemma we combine this last inequality with
(\ref{var0}). $\square$

\bigskip

Based on this lemma we shall provide another bound needed for the proof of
Theorem \ref{ThPeligrad}.

\begin{lemma}
\label{Lma bound var 2}Assume in addition to the conditions of Lemma
\ref{lmbound var} that the sequence $(\eta_{n})$ has the Markov property.
Then, for some universal constant $C$,
\begin{equation}
||\sup_{n}\frac{1}{n}E^{0}\left(  S_{n}^{2}\right)  ||_{1,w}\leq CE(V_{0}%
^{2})+C\sum\nolimits_{n\geq1}\frac{1}{n^{2}}E\left(  E(S_{n}|\eta_{0},\eta
_{n})\right)  ^{2}\text{.} \label{bound 2}%
\end{equation}

\end{lemma}

Proof. This bound follows from Lemma \ref{lmbound var}. We start by noting
that, by the properties of conditional expectations and the Markov property,
\begin{align*}
E\left(  S_{2^{k}}\bar{S}_{2^{k}}|\eta_{0}\right)   &  =E\left(  S_{2^{k}%
}E(\bar{S}_{2^{k}}|\eta_{2^{k}})|\eta_{0}\right)  =E\left(  E(S_{2^{k}}%
E(\bar{S}_{2^{k}}|\eta_{2^{k}})|\eta_{0},\eta_{2^{k}})|\eta_{0}\right) \\
&  =E\left(  E(S_{2^{k}}|\eta_{0},\eta_{2^{k}})E(\bar{S}_{2^{k}}|\eta_{2^{k}%
})|\eta_{0}\right)  .
\end{align*}
So, by the Cauchy-Schwartz inequality,%
\begin{align*}
E|E\left(  S_{2^{k}}\bar{S}_{2^{k}}|\eta_{0}\right)  |  &  \leq E|E(S_{2^{k}%
}|\eta_{0},\eta_{2^{k}})E(\bar{S}_{2_{k}}|\eta_{2^{k}})|\\
&  \leq\frac{1}{2}E\left(  E(S_{2^{k}}|\eta_{0},\eta_{2^{k}})\right)
^{2}+\frac{1}{2}E\left(  E(\bar{S}_{2_{k}}|\eta_{2^{k}})\right)  ^{2}\\
&  \leq E\left(  E(S_{2^{k}}|\eta_{0},\eta_{2^{k}})\right)  ^{2}.
\end{align*}
Therefore%
\[
\sum\nolimits_{k\geq0}\frac{1}{2^{k}}E|E\left(  (S_{2^{k}}\bar{S}_{2^{k}%
})|\eta_{0}\right)  |\leq\sum\nolimits_{k\geq0}\frac{1}{2^{k}}E\left(
E(S_{2^{k}}|\eta_{0},\eta_{2^{k}})\right)  ^{2}.
\]
As proven in Lemmas 12 and 13 in \cite{Pel20b}, for some positive constant
$c$,
\[
\sum\nolimits_{k\geq0}\frac{1}{2^{k}}E\left(  E(S_{2^{k}}|\eta_{0},\eta
_{2^{k}})\right)  ^{2}\leq c\sum\nolimits_{n\geq1}\frac{1}{n^{2}}%
E(E(S_{n}|\eta_{0},\eta_{n}))^{2}.
\]
It remains to apply Lemma \ref{lmbound var} to obtain the desired result.
$\square$

\bigskip

\textbf{Proof of Theorem \ref{ThPeligrad}}

\bigskip

The CLT and the convergence of moments under condition (\ref{condpf}) are
known (see Corollary 9 in \cite{Pel20b}). The proof of the quenched CLT is
based on the proof of Proposition \ref{pr quenched} combined with Lemma
\ref{Lma bound var 2}.

For $m$ fixed, we apply Lemma \ref{Lma bound var 2} with $\eta_{\ell+1}%
=(\xi_{\ell m},\xi_{(\ell+1)m})$\ and the sequence $V_{\ell+1}(m)=E(Y_{\ell
}|\xi_{\ell m},\xi_{(\ell+1)m})/\sqrt{m}$ where $(Y_{\ell})_{\ell\in Z}$ is
the extension to a stationary sequence of $Y_{k}$ defined in (\ref{defY}). It
is easy to see that, by using the Markov property and the properties of the
conditional expectation, we obtain for $k\geq0$%
\[
E(\sum\nolimits_{j=1}^{k+1}V_{j}|\eta_{0},\eta_{k+1})=\frac{1}{\sqrt{m}%
}E(S_{km}|\xi_{0},\xi_{km})+V_{k+1}.
\]
It follows that%
\[
||E(\sum\nolimits_{j=1}^{k+1}V_{j}|\eta_{0},\eta_{k+1})||^{2}\leq\frac{2}%
{m}||E(S_{km}|\xi_{0},\xi_{km})||^{2}+\frac{2}{m}||E(S_{m}|\xi_{0},\xi
_{m})||^{2}.
\]
So, for $R_{u}(m)$ defined in (\ref{defRn}), $R_{u}(m)=\sum\nolimits_{j=1}%
^{u}V_{j}(m),$ we obtain by Lemma \ref{Lma bound var 2}, for some $C_{1}>0,$%

\[
||\sup_{u}\frac{1}{u}E^{0}(R_{u}^{2}(m))||_{1,w}\leq C_{1}\sum\nolimits_{k=1}%
^{\infty}\frac{1}{k^{2}m}E(E(S_{km}|\xi_{0},\xi_{km}))^{2}.
\]
By the Cauchy-Schwartz inequality, and the properties of the conditional
expectation,%
\[
\frac{1}{mk^{2}}E(E(S_{km}|\xi_{0},\xi_{km}))^{2}\leq\frac{1}{k^{2}}%
E(E(S_{k}|\xi_{0},\xi_{k}))^{2}%
\]
and also%
\[
\sum\nolimits_{k=1}^{\infty}\frac{1}{mk^{2}}E(E(S_{km}|\xi_{0},\xi_{km}%
))^{2}\leq\sum\nolimits_{k=1}^{\infty}\frac{1}{k^{2}}E(E(S_{k}|\xi_{0},\xi
_{k}))^{2}<\infty.
\]
For any $k$ fixed, by (\ref{condpf}), we have that
\[
\lim_{m\rightarrow\infty}\frac{1}{mk}E(E(S_{km}|\xi_{0},\xi_{km}))^{2}=0.
\]
So, by the dominated convergence theorem for discrete measures,%
\[
\sum\nolimits_{k=1}^{\infty}\frac{1}{mk^{2}}E(E(S_{km}|\xi_{0},\xi_{km}%
))^{2}\rightarrow0\text{ as }m\rightarrow\infty.
\]
It follows that%
\[
\lim_{m\rightarrow\infty}||\sup_{u}\frac{1}{u}E^{0}(R_{u}^{2}(m))||_{1,w}=0.
\]
By Theorem 4.1 in Billingsley \cite{Bil}, note that the Fatou Lemma also holds
in the space $L^{1,w}$.\ Therefore,%
\[
||\lim\inf_{m\rightarrow\infty}\sup_{u}\frac{1}{u}E^{0}(R_{u}^{2}%
(m))||_{1,w}\leq\lim_{m\rightarrow\infty}||\sup_{u}\frac{1}{u}E^{0}(R_{u}%
^{2}(m))||_{1,w}=0.
\]
and so%
\[
\lim\inf_{m\rightarrow\infty}\sup_{u}\frac{1}{u}E^{0}(R_{u}^{2}(m))=0\text{
a.s.}%
\]
This proves that the martingale decomposition\ in (\ref{negli1}) holds. The
proof is now ended as in the proof of Proposition \ref{pr quenched}. $\square$

\bigskip

\begin{acknowledgement}
This paper was partially supported by the NSF grant DMS-2054598. It is
dedicated to Michael Lin's 80th birthday.
\end{acknowledgement}

\end{document}